\documentclass[11pt,a4paper]{amsart}
\usepackage[T1]{fontenc}
\usepackage[utf8]{inputenc}
\usepackage{lmodern,amssymb,amscd,mathtools}
\usepackage[margin=28mm]{geometry}
\usepackage{booktabs,array,enumitem,microtype}
\usepackage{tikz}
\usepackage[colorlinks=true,linkcolor=blue,citecolor=blue,urlcolor=blue]{hyperref}
\usepackage{bookmark}
\hypersetup{pdftitle={A compact nonformal manifold with holonomy Spin(7)},
 pdfauthor={Yaroslav Bazaikin},pdfsubject={Exceptional holonomy and Massey products}}
\numberwithin{equation}{section}
\newtheorem{theorem}{Theorem}[section]
\newtheorem{proposition}[theorem]{Proposition}
\newtheorem{lemma}[theorem]{Lemma}
\newtheorem{corollary}[theorem]{Corollary}
\theoremstyle{definition}
\newtheorem{definition}[theorem]{Definition}
\theoremstyle{remark}

\newcommand{\R}{\mathbb R}
\newcommand{\Q}{\mathbb Q}
\newcommand{\Z}{\mathbb Z}
\newcommand{\C}{\mathbb C}
\newcommand{\Spin}{\operatorname{Spin}}

\newcommand{\EH}{\operatorname{EH}}
\newcommand{\id}{\operatorname{id}}
\newcommand{\Fix}{\operatorname{Fix}}
\newcommand{\PD}{\operatorname{PD}}

\newcommand{\sgn}{\operatorname{sign}}
\newcommand{\Apl}{A_{\rm PL}}
\newcommand{\Omm}{\mathcal O}
\newcommand{\Mmm}{M}
\newcommand{\CP}{\mathbb{CP}}
\newcommand{\RP}{\mathbb{RP}}
\newcommand{\vol}{\operatorname{vol}}
\newcommand{\MM}{Mart\'in-Merch\'an}
\setlist{itemsep=2pt,topsep=4pt}

\title[Nonformality and Spin(7) holonomy]
 {A compact nonformal manifold with holonomy $\Spin(7)$}
\author{Yaroslav Bazaikin}
\address{Department of Mathematics, Faculty of Science,
  Jan Evangelista Purkyně University in Ústí nad Labem, Czech Republic}
\date{15 November 2026}
\makeatletter
\@ifundefined{subjclassname@2020}{%
  \@namedef{subjclassname@2020}{\textup{2020} Mathematics Subject Classification}}{}
\makeatother
\begin{document}
\begin{abstract}
We construct a compact simply connected eight-manifold with holonomy
$\Spin(7)$ and a triple Massey product of degrees $(3,2,2)$ that does
not contain zero.
The manifold is a resolution of a flat orbifold. Joyce's gluing theorem
gives it a torsion-free $\Spin(7)$ structure. The resolution map
factors through a quotient of $Y\times S^1$, where $Y$ is
\MM's nonformal $G_2$ manifold. The map to this quotient has degree
one, so pullback preserves a nonzero Massey pairing. On the resolution
the detecting degree-two class annihilates the full indeterminacy.
Simple connectivity and nonformality imply full $\Spin(7)$ holonomy.
The manifold has $(b_2,b_3,b_4^+,b_4^-)=(22,10,91,27)$.
\end{abstract}

\maketitle

\section{Introduction}\label{sec:status}
Compact K\"ahler manifolds are formal by the theorem of Deligne,
Griffiths, Morgan and Sullivan~\cite[Section~6]{DGMS}.
Cavalcanti~\cite{Cavalcanti} studied topological conditions for formality
in the context of special holonomy. Amann and Kapovitch recorded the
conjecture that every compact simply connected manifold with special
holonomy is formal~\cite[Conjecture~2]{AK}.
For holonomy $G_2$, Amann and Taimanov proved formality for one of
Joyce's examples~\cite{AT}, and \MM\ constructed a compact simply
connected nonformal example~\cite{MM}. Fiorenza, Kawai, L\^e and Schwachh\"ofer
showed that on a compact torsion-free $\Spin(7)$-manifold, a nontrivial
real triple Massey product must have degrees $(2,2,2)$ or a permutation
of $(2,2,3)$~\cite[Corollary~5.5(3) and Remark~5.6]{FKLS}. We realize the case
$(3,2,2)$ on a compact simply connected manifold with full holonomy
$\Spin(7)$.

\begin{theorem}\label{thm:main}
There exists a compact simply connected smooth eight-manifold $\Mmm$
admitting a torsion-free $\Spin(7)$ structure whose metric has holonomy
equal to $\Spin(7)$. There are rational classes
$A\in H^3(\Mmm)$ and $B,C,W\in H^2(\Mmm)$ such that
$\langle A,B,C\rangle$ is defined and
\begin{equation}\label{eq:mainperiod}
 \left|\left\langle nW,[\Mmm]\right\rangle\right|=8
 \qquad\text{for every }n\in\langle A,B,C\rangle.
\end{equation}
In particular, $\Mmm$ is nonformal over both $\Q$ and $\R$.
\end{theorem}

Let $\rho:Y\to X$ be the resolution of the flat seven-dimensional
orbifold in \MM's example. We choose an involution $\tau$ of $X$
that lifts to an orientation-reversing involution $\widetilde\tau$
of $Y$. Combining $\tau$ with the circle reflection $s\mapsto-s$
gives the flat orbifold $\Omm=(X\times S^1)/\langle(\tau,-1)\rangle$.
We construct a smooth resolution of $\Omm$ whose map factors as
\[
 \Mmm\xrightarrow{\ r\ } Z\xrightarrow{\ \bar\rho\ }\Omm,\qquad
 Z=(Y\times S^1)/\langle(\widetilde\tau,-1)\rangle .
\]
The map $\bar\rho$ is induced by $\rho\times\id$ and resolves the
branches inherited from $X$. The map $r$ resolves the remaining
singularities of $Z$, which lie along disjoint four-dimensional centers.

For the metric we use the flat presentation $\Omm=T^8/\Gamma$.
We construct a complex resolution $V\to T^8/H$ for an index-two
subgroup $H\subset\Gamma$, resolving the inherited branch first at
each intersection. The residual quotient of $V$ has a single real
singular component, whose resolution gives $\Mmm$. These local models
carry compatible quasi-asymptotically locally Euclidean (QALE)
$\Spin(7)$ structures, so Joyce's gluing theorem gives a torsion-free
$\Spin(7)$ structure on $\Mmm$.

The first Massey input on $Y$ and the circle form $ds$ both change
sign under the product involution, so their product is invariant.
A defining system on $Y\times S^1$ descends to $Z$ and pulls back
along the degree-one map $r$ to $\Mmm$, where a degree-two class
detects the product and annihilates its full indeterminacy.
The quotient $Z$ is simply connected, and
resolving its remaining centers preserves its fundamental group.
Nonformality then rules out the remaining proper holonomy groups,
which would make $\Mmm$ K\"ahler.

\begin{corollary}\label{cor:numerical}
The manifold in Theorem~\ref{thm:main} has
\begin{align}
 (b_k)_{k=0}^8&=(1,0,22,10,118,10,22,0,1),\label{eq:mainbetti}\\
 (b_2,b_3,b_4^+,b_4^-)&=(22,10,91,27),\nonumber\\
 \chi(\Mmm)&=144,\qquad\sgn(\Mmm)=64.\label{eq:mainindex}
\end{align}
\end{corollary}

After I submitted the first version of this paper, Jonas Stelzig pointed
me to a September 2026 note posted by Aleksandar Milivojevi\'c~\cite{MilivojevicSpin7}.
It describes a construction of a compact simply connected nonformal
manifold with full holonomy $\Spin(7)$ using the same affine involution
$\tau$ and hence the same flat orbifold $\Omm$. Both constructions
resolve the inherited branch first where three branches meet and
multiply the first entry of \MM's Massey product by the circle class
before passing to the quotient. In the present paper, we use the global
complex resolution $V\to T^8/H$ to construct the metric. We also
derive~\eqref{eq:mainperiod} from the rational transfer formula in
Proposition~\ref{prop:abstracttransfer} and compute the Betti numbers
and signature of $\Mmm$.

Sections~\ref{sec:affine} and~\ref{sec:rdata} describe the flat quotient
and construct its torsion-free resolution.
Section~\ref{sec:comparison} proves the factorization and determines
the action of $\widetilde\tau$ on the exceptional classes.
Section~\ref{sec:topology} completes the proof of
Theorem~\ref{thm:main}. Appendix~\ref{sec:mm} fixes the sign and
normalization of the seven-dimensional Massey period.
Appendix~\ref{sec:manualorbits} computes the Betti numbers and signature.

All cohomology groups have rational coefficients unless stated otherwise.
Products of differential forms denote wedge products. The orientations
are $dt\,dx_1\,dy_1\,dx_2\,dy_2\,dx_3\,dy_3$ in dimension seven and
$ds\wedge\vol_7$ on its product with the last circle.

\section{Massey products under an involution}\label{sec:foundations}
A space is formal over a field of characteristic zero if its algebra of
polynomial differential forms is connected by multiplicative
quasi-isomorphisms to its cohomology algebra with zero differential.
For a smooth manifold over $\R$, one may use the de Rham algebra
instead~\cite{Sullivan}. A defined triple Massey product which does not
contain zero is an obstruction to formality.

We use the following convention. If closed elements $\alpha,\beta,\gamma$
have degrees $p,q,r$ and $du=\alpha\beta$, $dv=\beta\gamma$, then
\begin{equation}\label{eq:masseyconvention}
 m=\alpha v+(-1)^{p+1}u\gamma
\end{equation}
is closed and represents the triple product. Its indeterminacy is
\begin{equation}\label{eq:generalindet}
 [\alpha]H^{q+r-1}+H^{p+q-1}[\gamma].
\end{equation}

\begin{lemma}\label{lem:detector}
Let $P$ be a closed oriented rational homology $n$-manifold. If
$w\in H^{n-p-q-r+1}(P)$ satisfies $[\alpha]w=[\gamma]w=0$,
then pairing with $w$ is constant on $\langle[\alpha],[\beta],[\gamma]\rangle$.
A nonzero value implies that this product does not contain zero.
\end{lemma}
\begin{proof}
For arbitrary classes $h,k$ of the degrees in~\eqref{eq:generalindet},
associativity and graded commutativity give $(\alpha h+k\gamma)w=0$.
\end{proof}

For a finite group $K$ acting smoothly on a manifold $P$, choose
an equivariant triangulation~\cite{Illman} and take barycentric
subdivisions until the action is regular~\cite[Chapter~III, Section~1]{Bredon}.
The quotient triangulates $P/K$, and invariant rational piecewise
polynomial forms descend simplex by simplex. This identifies the
algebra of forms on the quotient with the invariant subalgebra.
Here $\Apl(X)$ denotes the algebra of rational polynomial forms on the
singular simplices of $X$. Restricting to the simplices of a chosen
triangulation gives a multiplicative quasi-isomorphism to the algebra of
piecewise polynomial forms on that triangulation.
These algebras therefore compute rational cohomology and
Massey products~\cite[Section~7]{Sullivan}. Averaging over $K$ then gives
$H^*(P/K;\Q)=H^*(P;\Q)^K$.
If the action preserves orientation, its local links are finite
quotients of spheres with the rational cohomology of spheres.
The quotient is therefore a rational homology manifold. For an
effective action with trivial generic stabilizer, integration satisfies
\begin{equation}\label{eq:quotientintegration}
 \int_{P/K}\eta=\frac1{|K|}\int_P\widetilde\eta,
\end{equation}
where $\eta$ is a compactly supported top-degree form on $P/K$ and
$\widetilde\eta$ is its pullback to $P$.
For continuous maps we use the naturality of Massey products via
the functor $\Apl$.

\begin{proposition}\label{prop:abstracttransfer}
Let $Y$ be a closed oriented seven-manifold with an orientation-reversing
smooth involution $\tau$. Suppose $a,b,c,w\in H^2(Y;\Q)$ satisfy
\[
 \tau^*(a,b,c,w)=(-a,b,c,w),\qquad ab=bc=aw=cw=0,
\]
and let
$\mu=\langle mw,[Y]\rangle\ne0$ for $m\in\langle a,b,c\rangle$.
Set $Z=(Y\times S^1)/\langle(\tau,-1)\rangle$ and orient the product by
$ds\wedge\vol_Y$, with $\int_{S^1}ds=1$.
If $r:M\to Z$ is a continuous map from a closed oriented eight-manifold
of degree $d\ne0$, there are rational classes
$A\in H^3(M)$ and $B,C,W\in H^2(M)$ for which
\begin{equation}\label{eq:transferperiod}
 \left\langle nW,[M]\right\rangle=\frac d2\,\mu
 \qquad(n\in\langle A,B,C\rangle).
\end{equation}
In particular, $M$ is nonformal over $\Q$ and $\R$.
\end{proposition}
\begin{proof}
We use rational piecewise polynomial forms on a $\tau$-equivariant
triangulation of $Y$ and the affine subdivision of $S^1=\R/\Z$
with vertices at $s=0,\tfrac14,\tfrac12,\tfrac34$.
For the product, take iterated barycentric subdivisions of the
product cell complex until the action of $(\tau,-1)$ is regular.

Choose closed representatives of $a,b,c,w$ with the prescribed
$\tau$ parities and primitives $du=ab$, $dv=bc$. Replacing $u,v$
by $(u-\tau^*u)/2$, $(v+\tau^*v)/2$, respectively, makes $u$ odd
and $v$ even without changing their differentials. The resulting
representative $m=av-uc$ is odd, so $ds\wedge m$ descends to $Z$.
By Lemma~\ref{lem:detector}, every class in $\langle a,b,c\rangle$,
including $[m]$, pairs with $w$ to give $\mu$.

Put $\theta=ds$ and set on $Y\times S^1$
\[
 A_0=a\theta,\qquad B_0=b,\qquad C_0=c,\qquad W_0=w.
\]
Primitives for $A_0B_0$ and $B_0C_0$ are $u\theta$ and $v$,
giving the Massey representative
\[
 n_0=A_0v+(u\theta)C_0=\theta(av-uc)=\theta m.
\]
The inputs and primitives are invariant under $\sigma=(\tau,-1)$,
so the entire defining system descends to $Z$. We retain the
notation for the descended forms. Let $q:Y\times S^1\to Z$ be
the quotient map. With the prescribed orientation and
$\int_{S^1}\theta=1$, equation~\eqref{eq:quotientintegration} gives
\[
 \langle n_0W_0,[Z]\rangle
 =\tfrac12\langle\theta mw,[Y\times S^1]\rangle=\tfrac\mu2.
\]

Set $(A,B,C,W)=r^*([A_0],[B_0],[C_0],[W_0])$.
By naturality, $n=r^*[n_0]\in\langle A,B,C\rangle$, and the degree
formula gives $\langle nW,[M]\rangle=d\mu/2$.
Since $q^*$ is injective in rational cohomology, the relations
$aw=cw=0$ imply $[A_0][W_0]=[C_0][W_0]=0$ on $Z$.
Thus $AW=CW=0$, and $W$ annihilates the full indeterminacy
\begin{equation}\label{eq:fullindet}
 A H^3(M)+H^4(M)C.
\end{equation}
Every rational representative therefore has the same nonzero pairing.
After extension to $\R$, the same representative and vanishing
relations give the same conclusion for every real representative.
Hence $M$ is nonformal over both fields.
\end{proof}

\section{The flat orbifold and its singularities}\label{sec:affine}
We identify the fixed-set configuration and normal quotient models
needed for Joyce's resolution construction.

\subsection{Tori, lattices and forms}
Let
\begin{equation}\label{eq:tori}
 E=\C/(\Z+i\Z),\qquad T^7=(\R/4\Z)_t\times E^3,
 \qquad T^8=T^7\times(\R/\Z)_s,
\end{equation}
where $z_j=x_j+iy_j$. The metric is the Euclidean metric in the unscaled real coordinates $(t,x_1,y_1,x_2,y_2,x_3,y_3,s)$. Set
\begin{align}
 \omega&=dx_1dy_1+dx_2dy_2+dx_3dy_3,\nonumber\\
 \varphi&=dt\wedge\omega+\Re(dz_1dz_2dz_3),\label{eq:phi}\\
 \Phi&=ds\wedge\varphi+*_7\varphi.\label{eq:Phi}
\end{align}
In particular,
\[
\varphi=dt(dx_1dy_1+dx_2dy_2+dx_3dy_3)
 +dx_1dx_2dx_3-dx_1dy_2dy_3-dy_1dx_2dy_3-dy_1dy_2dx_3.
\]
For the groups defined below, $\varphi$ descends to $T^7/G$ and $\Phi$ to $T^8/\Gamma$.

\subsection{The group action}
Define
\begin{align}
 F(t,z_1,z_2,z_3)&=(t+1,iz_1,iz_2,-z_3),\nonumber\\
 \iota_1(t,z_1,z_2,z_3)&=(t,-z_1,-z_2,z_3+\tfrac12),\nonumber\\
 \iota_2(t,z_1,z_2,z_3)&=(t,-z_1,-z_2,z_3+\tfrac i2),\label{eq:oldactions}\\
 \kappa(t,z_1,z_2,z_3)&=(1-t,\bar z_1,\bar z_2,\bar z_3).\nonumber
\end{align}
They act trivially on $s$. The group $G=\langle F,\kappa,\iota_1,\iota_2\rangle$ has order $32$ and preserves $\varphi$. The quotient $X=T^7/G$ is the orbifold of \MM: quotienting first by $F$ gives the mapping torus used there, then by $\iota_1,\iota_2$, and finally by $\kappa$.

Now put
\begin{equation}\label{eq:tau}
 \tau(t,z_1,z_2,z_3)=(1-t,z_1,z_2,\tfrac14-z_3),
 \qquad \sigma(t,z,s)=(\tau(t,z),-s).
\end{equation}
All translations are taken in the lattices~\eqref{eq:tori}. Direct composition gives
\begin{equation}\label{eq:relations}
 \tau^2=1,\quad \tau F\tau=F^{-1}\iota_1,\quad
 [\tau,\kappa]=[\tau,\iota_1]=[\tau,\iota_2]=1.
\end{equation}
Moreover $\tau^*\varphi=-\varphi$. This can be seen directly: $dt$ changes sign; the first two complex coordinates are unchanged, whereas $dz_3$ changes sign; and $dx_3dy_3$ is unchanged. The induced metric is preserved and the seven-dimensional orientation is reversed. It follows that $\tau^*(*_7\varphi)=*_7\varphi$, since reversal of orientation changes the Hodge star by a minus sign. Hence $\sigma^*\Phi=\Phi$.

\begin{definition}\label{def:orbifold}
The orbifold studied here is
\[
 \Omm=T^8/\Gamma,\qquad \Gamma=\langle G,\sigma\rangle.
\]
A singular branch is called \emph{old} if its order-two stabilizer belongs to $G$ and \emph{new} otherwise. This marking is part of the subsequent resolution prescription.
\end{definition}
Every group element has a unique normal form
\begin{equation}\label{eq:normalform}
 F^a\kappa^b\iota_1^c\iota_2^d\sigma^e,
 \qquad 0\leq a<4,\quad b,c,d,e\in\{0,1\}.
\end{equation}
For existence, move $\sigma$ to the right using~\eqref{eq:relations}, and use $\kappa F\kappa=F^{-1}$ and centrality of $\iota_1,\iota_2$ in $G$. For uniqueness, a normal-form identity must have $e=0$ from its action on $s$, $b=0$ from the sign of $t$, then $a=0$ from the translation of $t$ modulo $4$, and $c=d=0$ from the translations of $z_3$. Thus $|\Gamma|=64$, $G$ is normal of index two, and the old/new marking is conjugation invariant.
Since $T^8/G=X\times S^1$, taking the quotient by $\Gamma/G$ gives
$\Omm=(X\times S^1)/\langle(\tau,-1)\rangle$.

The subgroup
\begin{equation}\label{eq:H}
 H=\langle F^2,\kappa,\iota_1,\iota_2,\sigma\rangle
\end{equation}
consists of the elements with even $F$ exponent and has index two.
Its five displayed generators commute and have order two; the normal
form identifies it with $(\Z/2)^5$.
Arithmetic in binary indices and labels is over $\mathbb F_2=\Z/2\Z$.
In exponents, affine coordinates and integer-valued parameters, we use
their representatives $0,1$.
For $\alpha,c,d\in\mathbb F_2$, put
\begin{equation}\label{eq:OUJ}
\begin{aligned}
 O_{\alpha d}&=F^{2\alpha}\kappa\iota_2^d,
 &U_{\alpha cd}&=F^{2\alpha}\iota_1^c\iota_2^d\sigma,\\
 R_{\alpha c}&=F^{2\alpha+1}\kappa\iota_1^c,
 &J_c&=\kappa\iota_1^c\sigma.
\end{aligned}
\end{equation}

\subsection{The fixed-set configuration}\label{sec:fixedequations}
Following Joyce~\cite[Definition~13.2.1, p.~346]{JoyceBook}, let
\[
 \mathcal L=\{D:D\text{ is a connected component of }\Fix(A),
                 \ A\leq\Gamma\}.
\]
This includes $T^8=\Fix(\{1\})$. For $D\in\mathcal L$, its
pointwise stabilizer $A_D=\{\gamma\in\Gamma:\gamma|_D=\id\}$
acts on its constant normal space $W_D$. At a generic point of $D$
the normal quotient is $W_D/A_D$. The setwise stabilizer of $D$, modulo $A_D$, acts on
$D$ and on this quotient; the normal action is given by derivatives
of the affine maps. We determine $\mathcal L$ by listing the fixed
sets of individual elements and then their intersections.

Write $q=t/4\in\R/\Z$ for the affine coordinates; this does not change the metric convention. Decompose the coordinates into
\[
 \mathsf A=(x_1,y_1,x_2,y_2),\qquad \mathsf B=(q,x_3,y_3,s).
\]
For the element~\eqref{eq:normalform}, its action is
\begin{align}
 z_j&\longmapsto i^a(-1)^{c+d}\begin{cases}z_j&b=0,\\\bar z_j&b=1,\end{cases}
 &&(j=1,2),\label{eq:Aaction}\\
 q&\longmapsto(-1)^{b+e}\bigl(q-\tfrac18\bigr)+\tfrac18+\tfrac a4,\nonumber\\
 x_3&\longmapsto(-1)^{a+e}x_3+(-1)^a(\tfrac e4+\tfrac c2),\label{eq:Baction}\\
 y_3&\longmapsto(-1)^{a+b+e}y_3+(-1)^{a+b}\tfrac d2,\nonumber\\
 s&\longmapsto(-1)^e s.\nonumber
\end{align}
Fixed sets are products of their $\mathsf A$ and $\mathsf B$ parts. A circle reflection has two fixed points; a nontrivial translation has none. On $E$, $z\mapsto\pm\bar z$ fixes two circles for either sign, whereas $z\mapsto i\bar z$ and $z\mapsto-i\bar z$ fix the single circles $x=y$ and $x=-y$, respectively.

\begin{lemma}[All nonempty fixed sets]\label{lem:fixedelements}
Among the $63$ nonidentity elements of $\Gamma$, exactly $18$ have fixed points: $14$ have four-dimensional fixed tori and four have isolated fixed points. All $18$ are involutions. Their fixed sets are listed in Table~\ref{tab:elements}; all other elements act without fixed points.
\end{lemma}
\begin{proof}
If $e=b=0$, the $\mathsf B$ equations force $a=c=d=0$, giving the identity.

For $e=0,b=1$, they force $c=0$ for even $a$, and $d=0$ for odd $a$. In both cases, $q$ and one coordinate of $z_3$ are reflected, so the fixed set in $\mathsf B$ consists of four two-tori. Each $E$ factor in $\mathsf A$ contributes two circles for even $a$ and one for odd $a$, giving $16$ and $4$ four-tori, respectively.

For $e=1,b=0$, odd $a$ is excluded by a nonzero quarter translation in $x_3$. For $a=2\alpha$, all four $\mathsf B$ coordinates are reflected, while the action on $\mathsf A$ is $(-1)^{\alpha+c+d}$. Thus the fixed set consists of $16$ four-tori or $256$ points according as $\alpha+c+d=0$ or $1$.

Finally, $e=b=1$ forces $a=d=0$. The coordinates $x_3,s$ are reflected and $q,y_3$ are free; each $E$ factor contributes two fixed circles. This gives $16$ four-tori for each $J_c$. We obtain Table~\ref{tab:elements}, and all listed maps are involutions.
\end{proof}

\begin{table}[ht]
\centering
\caption{Nonidentity elements with fixed points. Components are counted
on $T^8$ for each element.}\label{tab:elements}
\begin{tabular}{@{}lcccl@{}}\toprule
Family & Elements & $\dim_{\R}$ & Components & Branch\\\midrule
$O_{\alpha d}$ & $4$ & $4$ & $16$ & old, complex\\
$R_{\alpha c}$ & $4$ & $4$ & $4$ & old, real\\
$U_{\alpha cd}$, $\alpha+c+d=0$ & $4$ & $4$ & $16$ & new\\
$U_{\alpha cd}$, $\alpha+c+d=1$ & $4$ & $0$ & $256$ & new\\
$J_c$ & $2$ & $4$ & $16$ & new\\\bottomrule
\end{tabular}
\end{table}

\subsection{Intersections and normal models}\label{sec:stabilizers}
Distinct old fixed tori are disjoint. This follows directly from~\eqref{eq:Baction}: different even $a$ have disjoint $t$ levels, different $d$ have disjoint $y_3$ levels; odd and even $a$ have disjoint $t$ levels; for fixed odd $a$, the two choices of $c$ have disjoint $x_3$ levels. Thus every nontrivial old stabilizer is cyclic of order two.

\begin{lemma}[Stabilizer bound]\label{lem:orderbound}
Every point stabilizer in $\Gamma$ has order $1$, $2$ or $4$. A stabilizer of order four is a Klein four group with exactly one nonidentity element in $G$.
\end{lemma}
\begin{proof}
For a stabilizer $K$, the intersection $K\cap G$ has order at most two. The quotient $K/(K\cap G)$ embeds in $\Gamma/G$, which has order two. Hence $|K|\leq4$. Every nonidentity element with a fixed point is an involution by Lemma~\ref{lem:fixedelements}, excluding a cyclic group of order four. If $|K|=4$, its intersection with $G$ must have order two.
\end{proof}
Any subgroup fixing $x$ lies in its point stabilizer $\Gamma_x$.
Thus, to complete $\mathcal L$, it remains to find the common
fixed sets of Klein four groups.
The old odd-$a$ fixed sets cannot intersect any new fixed set: their $x_3$ coordinates are multiples of $1/4$, whereas every new fixed point has $x_3$ an odd multiple of $1/8$. Thus every order-four stabilizer lies in the even-$a$ subgroup.

Here is an explicit formula for all eight order-four stabilizers. Choose $a\in\{0,2\}$ and $c,d\in\{0,1\}$, and put
\begin{equation}\label{eq:Kformula}
 g=F^a\kappa\iota_2^d,\qquad
 h=F^a\iota_1^c\iota_2^d\sigma,
 \qquad K(a,c,d)=\{1,g,h,gh\}.
\end{equation}
Writing $a=2\alpha$, we have $g=O_{\alpha d}$ and $h=U_{\alpha cd}$.
The $q$ and $y_3$ equations for a common fixed point force the
same $a$ and $d$, respectively, in $g$ and $h$. Conversely these
choices give nonempty common fixed loci. Since $gh=J_c$,
replacing $h$ by $J_c$ leaves the common fixed set and the
generated group unchanged.

Each of the four $\mathsf B$-coordinates has two possible fixed values, giving $16$ choices in total. On $\mathsf A$, $h$ acts by $(-1)^{a/2+c+d}$. If this sign is $+1$, the two real fixed directions of $g$ remain free and the common fixed set is a union of $64$ two-tori. If it is $-1$, all four $\mathsf A$ coordinates are fixed and there are $256$ points. This gives four groups of each type, hence $256$ two-tori and $1024$ points before taking the quotient.

For $\alpha+c+d=1$,
$\Fix(U_{\alpha cd})=\Fix(K(2\alpha,c,d))$;
hence every isolated fixed point has full stabilizer $K(2\alpha,c,d)$.

Write $A_1=\C^2/\{\pm1\}$. At an open four-dimensional stratum,
the stabilizer acts by $-\id$ on the normal space, so the full
local model is $\R^4\times A_1$.

At a two-dimensional fixed locus, split its normal space into the three complex character lines of $K$. The elements $g,h,gh$ act by
\begin{equation}\label{eq:curvecharacters}
 (-1,-1,1),\quad(-1,1,-1),\quad(1,-1,-1)
\end{equation}
on these lines, respectively. The full local model is $\R^2\times\C^3/(\Z/2)^2$.
Its three four-dimensional branches are fixed by $g,h,gh$, with
$g$ marked old.

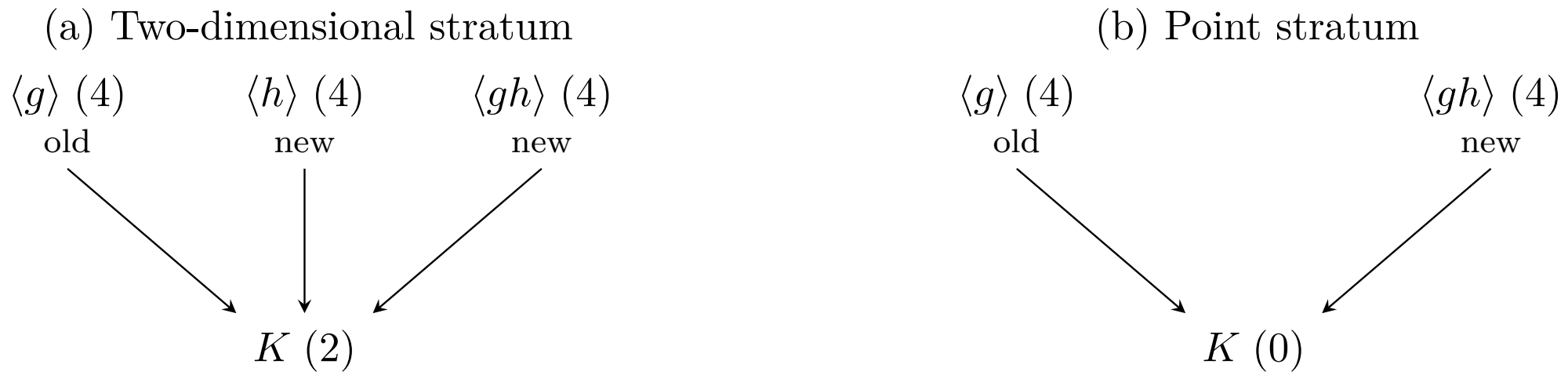
\begin{figure}[htbp]
\centering
\begin{minipage}[t]{.48\linewidth}
\centering\small
(a) Two-dimensional stratum\\[3pt]
\begin{tikzpicture}[font=\small,align=center,inner sep=2pt,
 outer sep=2pt,line width=0.45pt,>=stealth]
\node (g) at (-2,0) {$\langle g\rangle\;(4)$\\[-1pt]{\scriptsize old}};
\node (h) at (0,0) {$\langle h\rangle\;(4)$\\[-1pt]{\scriptsize new}};
\node (gh) at (2,0) {$\langle gh\rangle\;(4)$\\[-1pt]{\scriptsize new}};
\node (k) at (0,-2) {$K\;(2)$};
\draw[->] (g.south) -- (k.north west);
\draw[->] (h.south) -- (k.north);
\draw[->] (gh.south) -- (k.north east);
\end{tikzpicture}
\end{minipage}\hfill
\begin{minipage}[t]{.48\linewidth}
\centering\small
(b) Point stratum\\[3pt]
\begin{tikzpicture}[font=\small,align=center,inner sep=2pt,
 outer sep=2pt,line width=0.45pt,>=stealth]
\node (g) at (-2,0) {$\langle g\rangle\;(4)$\\[-1pt]{\scriptsize old}};
\node (gh) at (2,0) {$\langle gh\rangle\;(4)$\\[-1pt]{\scriptsize new}};
\node (k) at (0,-2) {$K\;(0)$};
\draw[->] (g.south) -- (k.north west);
\draw[->] (gh.south) -- (k.north east);
\end{tikzpicture}
\end{minipage}
\caption{Local incidence near the two types of strata with stabilizer
of order four. Vertices show full stabilizers, with real stratum dimensions
in parentheses; an arrow $S\to T$ means $T\subset\overline S$ locally.}
\label{fig:localstrata}
\end{figure}

At a point, $g$ acts as $(-1,-1,1,1)$ and $gh$ as $(1,1,-1,-1)$ in suitable complex coordinates, while $h$ acts as $-\id$ on $\C^4$. The local model is therefore
\begin{equation}\label{eq:pointmodel}
 (\C^2/\{\pm1\})\times(\C^2/\{\pm1\}).
\end{equation}
Here $h$ has an isolated fixed point, but the full stabilizer is $K$:
locally $\Fix(h)=\Fix(K)$. There are two four-dimensional branches,
fixed by $g$ and $gh$, marked old and new. Figure~\ref{fig:localstrata}
shows these incidences in $\mathcal L$.

The nontrivial members of $\mathcal L$ are the four-tori, two-tori
and points listed above. Distinct members with order-four pointwise
stabilizer are disjoint: an intersection
would force the stabilizers to coincide, and distinct components of
their fixed set cannot meet. Hence the only proper inclusions in the
singular part of $\mathcal L$ are those of two-tori or points in
four-tori. In particular, no point stratum lies on a two-dimensional
stratum. The $R_{\alpha c}$-loci are disjoint from the entire
configuration fixed by nontrivial subgroups of $H$.

The connected singular strata in $\Omm$ have dimensions $4,2,0$
and numbers $21,16,64$; the orbit counts are given in
Appendix~\ref{sec:manualorbits}. The affine action permutes
$\mathcal L$, preserves the old/new marking and determines the
normal transitions. We now construct compatible marked resolutions
and metrics for these data.

\section{The resolution and its R-data}\label{sec:rdata}
We use the index-two subgroup $H$ of~\eqref{eq:H} to construct the
resolution $\Mmm\to\Omm$. The intersecting branches of $T^8/H$ have
complex normal models. Resolving the old complex branches first and then
the new ones gives a smooth complex fourfold $V$. The residual action
of $F$ lifts to an involution $\bar F$ of $V$. Resolving the
separate real singular component of $V/\langle\bar F\rangle$ defines
$\Mmm$.

To obtain Joyce's R-data for the flat quotient $\Omm$, we supplement
the fixed-set data of Section~\ref{sec:affine} with local product
resolutions, torsion-free product $\Spin(7)$ structures with compatible
QALE limits, and lifts of the affine maps preserving these structures
\cite[Definition~13.2.2, p.~347]{JoyceBook}.

\subsection{A complex structure adapted to the Cayley form}
To resolve the complex models coherently, introduce the following
complex coordinates on the same real torus:
\begin{equation}\label{eq:ucoords}
 u_1=x_1+ix_2,\quad u_2=y_1-iy_2,\quad
 u_3=t-iy_3,\quad u_4=s+ix_3.
\end{equation}
Their lattices are, respectively,
\[
 \Z+i\Z,\quad\Z+i\Z,\quad4\Z+i\Z,\quad\Z+i\Z.
\]
Thus this is a product of four rectangular elliptic curves; in particular it is a projective complex torus. Its standard K\"ahler and volume forms are
\begin{equation}\label{eq:omegaprime}
 \omega'=dx_1dx_2-dy_1dy_2-dt\,dy_3+ds\,dx_3,
 \qquad\Omega'=du_1du_2du_3du_4.
\end{equation}
\begin{lemma}\label{lem:phase}
With the orientations fixed above,
\begin{equation}\label{eq:phase}
 \Phi=\tfrac12\omega'^2-\Re\Omega'.
\end{equation}
The subgroup $H$ of~\eqref{eq:H} is holomorphic and volume preserving. Its linear image is the group of all even sign changes of four complex coordinates. Every element outside $H$ is antiholomorphic.
\end{lemma}
\begin{proof}
Expanding~\eqref{eq:omegaprime} gives~\eqref{eq:phase}. In the $u$ coordinates the generators act as follows:
\begin{center}\small
\begin{tabular}{@{}cl@{}}\toprule
Map & Image of $(u_1,u_2,u_3,u_4)$\\\midrule
$F^2$ & $(-u_1,-u_2,u_3+2,u_4)$\\
$\kappa$ & $(u_1,-u_2,1-u_3,u_4)$\\
$\iota_1$ & $(-u_1,-u_2,u_3,u_4+i/2)$\\
$\iota_2$ & $(-u_1,-u_2,u_3-i/2,u_4)$\\
$\sigma$ & $(u_1,u_2,1-u_3,-u_4+i/4)$\\
$F$ & $(-\bar u_2,\bar u_1,\bar u_3+1,\bar u_4)$\\\bottomrule
\end{tabular}
\end{center}
The first five rows have complex determinant one. Their sign vectors span the even-sign subgroup of rank three; the kernel of the linear representation consists of translations. The last row conjugates the complex structure and sends $\omega'$ to $-\omega'$ and $\Omega'$ to $\bar\Omega'$, preserving~\eqref{eq:phase}.
\end{proof}
The quotient $T^8/H$ has only the complex singularities already classified. The residual map induced by $F$ is an antiholomorphic involution, since $F^2\in H$. Its fixed locus is disjoint from the complex singular locus: otherwise an odd-$F$ element and an element of $H$ would fix the same point, contrary to Section~\ref{sec:stabilizers}.

\subsection{Eguchi--Hanson replacements and markings}\label{sec:ehmarkings}
The minimal resolution of $A_1$ is
\[
 \rho_{\EH}:\operatorname{Tot}\mathcal O_{\CP^1}(-2)\longrightarrow A_1.
\]
It carries the Eguchi--Hanson ALE hyperk\"ahler metrics, and its
exceptional sphere $S$ has self-intersection $-2$. We write $\EH$ for
the underlying resolution. An oriented ray in the three-dimensional
space of self-dual normal two-forms determines a complex structure
$I$, normalized by $I^2=-1$, and orients $S$. The area of $S$ fixes
the scale of the metric. We call the generator $u\in H^2(\EH;\Z)$
with $\langle u,[S]\rangle=1$ the exceptional root. The signs with
which transitions act on $u$ define the root local system. If all
monodromy signs are positive, this system is trivial and the
exceptional cohomology summand is untwisted.

At a smooth complex $A_1$ stratum we use the complex structure inherited from~\eqref{eq:ucoords}. At an intersection, the marking by $G$ prescribes resolving the old branch first.

\subsection{Marked local resolutions}\label{sec:fan}
For the transverse curve model~\eqref{eq:curvecharacters}, the quotient lattice is
\begin{equation}\label{eq:fanlattice}
 N=\Z^3+\Z\tfrac{e_1+e_2}{2}+\Z\tfrac{e_1+e_3}{2},
 \qquad[N:\Z^3]=4.
\end{equation}
The invariant monomial lattice is
\[
 M_N=\{(r_1,r_2,r_3)\in\Z^3:r_1\equiv r_2\equiv r_3\pmod2\}.
\]
Set $m_{ij}=(e_i+e_j)/2$, with $m_{12}$ marked old. The four maximal cones are
\begin{equation}\label{eq:fancones}
 \begin{aligned}
 C_1&=\langle e_1,m_{12},m_{13}\rangle,&
 C_2&=\langle e_2,m_{12},m_{23}\rangle,\\
 C_3&=\langle e_3,m_{12},m_{13}\rangle,&
 C_4&=\langle e_3,m_{12},m_{23}\rangle.
 \end{aligned}
\end{equation}
They result from subdividing first along $m_{12}$, then along $m_{13}$ and $m_{23}$. Each determinant in $\Z^3$ has absolute value $1/4$, which is the covolume of $N$; hence each cone is regular in the quotient lattice. All primitive rays lie in $r_1+r_2+r_3=1$, so the resolution is crepant.

We construct a strictly convex support function $\psi$ by assigning heights
\begin{equation}\label{eq:heights}
 h(e_i)=0,\qquad h(m_{12})=-3,\qquad
 h(m_{13})=h(m_{23})=-1.
\end{equation}
For each cone $C_j$, let $\ell_j$ be the linear function matching these heights on its ray generators. These linear pieces agree on common faces and define $\psi$ by $\psi|_{C_j}=\ell_j$. The last column lists $h(v)-\ell_j(v)$ in the order $(e_1,e_2,e_3,m_{12},m_{13},m_{23})$:
\begin{center}\small
\begin{tabular}{@{}ccc@{}}\toprule
Cone & $\ell_j$ & Gaps\\\midrule
$C_1$ & $(0,-6,-2)$ & $(0,6,2,0,0,3)$\\
$C_2$ & $(-6,0,-2)$ & $(6,0,2,0,3,0)$\\
$C_3$ & $(-2,-4,0)$ & $(2,4,0,0,0,1)$\\
$C_4$ & $(-4,-2,0)$ & $(4,2,0,0,1,0)$\\\bottomrule
\end{tabular}
\end{center}
All off-cone gaps are positive, so the resolution is projective over the affine quotient. The fan and support function are invariant under interchange of the two new branches, equivalently of $e_1$ and $e_2$. The old/new marking is $\Gamma$-invariant, so the only possible nontrivial permutation of the three branches exchanges the two new ones.

The heights~\eqref{eq:heights} will be used only to establish
projectivity. We choose the metric classes separately in
Section~\ref{sec:QALE}.

\begin{lemma}[Resolving the old branch first]\label{lem:markedresolution}
In both maximal local models, resolving the old branch first leaves
smooth, disjoint new centers with transverse singularity $A_1$.
Resolving these centers gives the fan~\eqref{eq:fancones} in the curve
model and $\EH\times\EH$ in the point model. These resolutions and
the contractions of their new exceptional factors commute with
isomorphisms preserving the old branch, including antiholomorphic ones.
\end{lemma}
\begin{proof}
In the curve model, choose local transverse complex coordinates $(z_1,z_2,z_3)$ as in~\eqref{eq:curvecharacters}. The old and new generators then act by
\[
 g(z_1,z_2,z_3)=(-z_1,-z_2,z_3),\qquad
 h(z_1,z_2,z_3)=(-z_1,z_2,-z_3).
\]
After resolving $g$, the two charts of $\EH$ are
\[
 (v,w)=(z_1^2,z_2/z_1),\qquad
 (v',w')=(z_2^2,z_1/z_2),\qquad
 w'=w^{-1},\quad v'=vw^2.
\]
The residual involution acts by
\begin{equation}\label{eq:EHinvolution}
 (v,w,z_3)\longmapsto(v,-w,-z_3),\qquad
 (v',w',z_3)\longmapsto(v',-w',-z_3).
\end{equation}
Its fixed lines $w=z_3=0$ and $w'=z_3=0$ meet the exceptional sphere
at distinct poles and do not meet on the chart overlap.
Their images in the quotient have transverse singularity $A_1$.
Restoring the common tangent coordinate gives two disjoint smooth
surface centers. Resolving them inserts $m_{13}$ and $m_{23}$ after
$m_{12}$, producing~\eqref{eq:fancones} and the morphisms
\begin{equation}\label{eq:localcontractions}
 \widetilde{\C^3/K}\longrightarrow
 (\EH\times\C)/\langle h\rangle\longrightarrow\C^3/K.
\end{equation}

This first step also has an intrinsic description. The quotient is
\[
 \operatorname{Spec}\C[X,Y,Z,T]/(T^2-XYZ),\qquad
 X=z_1^2,\quad Y=z_2^2,\quad Z=z_3^2,\quad T=z_1z_2z_3.
\]
The old axis has ideal $(X,Y,T)$. On the positive cone its Newton
function is
\[
 \min(2n_1,2n_2,n_1+n_2+n_3)=\min(2n_1,2n_2).
\]
It changes linear expression across $n_1=n_2$, so the normalized
blowup of this ideal is the subdivision along $m_{12}$.
It is therefore the same partial resolution as
$(\EH\times\C)/\langle h\rangle$.
In a related surface construction~\cite[Section~3]{Bazaikin2005},
$M_{1,2}$ replaces a cyclic quotient singularity and $M_{1,1}$
resolves the single remaining $A_1$ point.

In the point model, $g=(-1,-1,1,1)$ and $gh=(1,1,-1,-1)$.
Resolving $g$ and quotienting by the residual action of $gh$
gives $\EH\times A_1$.
The remaining center is $\EH\times\{0\}$, and its minimal
resolution gives $\EH\times\EH$.

Normalized blowups of the old branch ideals commute with
isomorphisms preserving those ideals. The remaining minimal $A_1$
resolutions have the same property. The constructions commute with
complex conjugation as well. This proves the stated compatibility
of the resolutions and their contractions.
\end{proof}

\subsection{Global complex resolution}\label{sec:globalV}
We construct
\begin{equation}\label{eq:Vmap}
 f:V\longrightarrow T^8/H
\end{equation}
by resolving the old complex branches first and then the new ones.
The old branches are disjoint, since a full stabilizer contains at
most one nonidentity old element. Their normalized blowups can
therefore be performed simultaneously. At a generic old surface this
is the minimal $A_1$ resolution; at an intersection it is the first
step of Lemma~\ref{lem:markedresolution}. By that lemma the remaining
centers are smooth and disjoint, with normal singularity $A_1$.
Their minimal resolutions complete the construction.

The local models show that $V$ is smooth and crepant. It is also
projective: the base $T^8/H$ is projective and both steps are
projective blowups. The support function~\eqref{eq:heights} gives a
relatively ample exceptional divisor on the curve model; its
old and new coefficients agree with the corresponding choices on
the ordinary and product models.

On complex normal coordinates, the affine transitions of
Section~\ref{sec:affine} are signed complex permutations with unitary
phases, or their complex conjugates, preserving the old branch.
Diagonal phases act on the toric resolution, and the relevant
permutations preserve its fan. Lemma~\ref{lem:markedresolution}
therefore glues the local morphisms and lifts the involution induced
by $F$ on $T^8/H$ to an antiholomorphic map $\bar F$ on $V$.
The lifts satisfy the original group relations: two such maps
agreeing off the exceptional locus agree everywhere by continuity.
In particular, $\bar F$ is an involution.

By equivariance, $f$ descends to a partial resolution
$V/\langle\bar F\rangle\to\Omm$, since
$(T^8/H)/(\Gamma/H)=\Omm$. The remaining singular locus is the image
of $\Fix(\bar F)$; we describe its resolution next.

\subsection{The real resolution and the manifold \texorpdfstring{$\Mmm$}{M}}\label{sec:realL}
The involution $\bar F$ on $V$ has exactly one fixed component
\begin{equation}\label{eq:L}
 L\cong N_1\times S^1,\qquad
 N_1=T^3/\langle(a,b,c)\mapsto(-a,-b,c+\tfrac12)\rangle.
\end{equation}
Here $a,b,c$ have period one. The deck involution is free. Use the component $D$ of $\Fix(F^3\kappa)$ with $t=0$, $z_1=a-ia$, $z_2=b-ib$, $z_3=ic$. Write a group element as $F^j\kappa^k\iota_1^l\iota_2^m\sigma^n$. To preserve $D$, the $x_3$ equation forces $n=l=0$, and the $t$ equation forces $(j,k)=(0,0)$ or $(3,1)$; $m$ is arbitrary. All four resulting elements preserve $D$, so its setwise stabilizer is $\langle F^3\kappa,\iota_2\rangle$. Its orbit has $64/4=16$ components, which is the total number of old odd-$F$ components. This proves both connectedness of the residual locus and the asserted quotient description. The element $\iota_2$ acts as the displayed free involution and trivially on $s$. The tangent and normal bundles of $L$ have finite flat orthogonal holonomy.

No additional fixed point appears on the complex exceptional locus. If there were one, its image in $T^8/H$ would have an odd-$F$ stabilizer in addition to a nontrivial element of $H$, contradicting Section~\ref{sec:stabilizers}. Thus the complex blowups leave a neighborhood of $L$ unchanged, with its flat tangent and normal bundles.

In $V/\langle\bar F\rangle$, the image of $L$, still denoted by $L$, has transverse model $\R^4/\{\pm1\}$. Resolve this model using the parallel direction $\ell=dy_3=dc$ on $N_1$. With the induced product orientation, $ds\wedge\ell+*_3\ell$ fixes a self-dual direction on $L$, and the associated normal hyperk\"ahler direction specifies $\EH$. The element $\iota_2$ preserves that direction, and its unitary lift preserves the orientation of the exceptional sphere. Thus the root local system is trivial, giving the \emph{untwisted} choice.

The lift acts isometrically on the product model, so the difference between the parallel Cayley form and its pullback has constant norm. Since $\iota_2$ preserves the flat limit $\Phi$, this norm tends to zero along the Eguchi--Hanson end and hence vanishes. The size of the replacement along $L$ can be chosen independently, since $L$ is disjoint from the complex strata. We use the same size for the two old branches exchanged by $\tau$ in Section~\ref{sec:comparison}.

Let $\Mmm$ be the compact smooth manifold obtained from
$V/\langle\bar F\rangle$ by this untwisted Eguchi--Hanson replacement
along $L$. The gluing uses $\rho_{\EH}$ in the fixed flat tubular
coordinates. The resolution map $\pi_M\colon\Mmm\to\Omm$ is the composite
\[
 \Mmm\longrightarrow V/\langle\bar F\rangle\longrightarrow\Omm,
\]
where the second map is induced by $f$.

\subsection{Compatible metrics and Joyce's theorem}\label{sec:QALE}
It remains to equip the complex local resolutions defining $\Mmm$
with compatible metrics.
A QALE K\"ahler metric is asymptotic to the flat quotient away from its
singular directions and to products with lower-dimensional resolved
models along those directions~\cite[Definition~2.1]{QALE}.
Joyce proves that a crepant resolution
is a local product resolution and that each K\"ahler class admitting
a QALE representative contains a unique Ricci-flat QALE K\"ahler
metric~\cite[Proposition~3.2 and Theorem~3.3]{QALE}.
We now choose such classes with matching restrictions at all
incidences of strata.

For the curve model, Example~3.6 of~\cite{QALE} gives QALE
K\"ahler metrics on all four crepant resolutions of
$\C^3/(\Z/2)^2$, including the marked fan~\eqref{eq:fancones}.
The transitions induced by $\Gamma$ preserve the $U$ and $J$ families
separately by~\eqref{eq:Fconjugation}. Average one such K\"ahler form
over the finite group generated by the linear transition maps and the
symmetry $e_1\leftrightarrow e_2$ of the marked fan. Use $\gamma^*\omega$
for holomorphic maps and $-\gamma^*\omega$ for antiholomorphic maps.
Each summand is a positive closed $(1,1)$-form. Their flat limits
agree, and averaging their lower-dimensional limits preserves the
QALE conditions. Thus the averaged form is again QALE and has the
required symmetries.

In each asymptotic $A_1$ factor, the period of the averaged form over
the exceptional sphere is its area. Write $\lambda_{\rm old}>0$ for
the old period and $\lambda_{\rm new}>0$ for the common period of the
two new spheres exchanged by this fan symmetry. Use this same marked class
on every curve model. On a point model use the product of the
Eguchi--Hanson classes with periods
$\lambda_{\rm old}$ and $\lambda_{\rm new}$; on an ordinary complex
$A_1$ branch use the class with the corresponding period.
This is a single choice of periods across the strata.
A common dilation supplies the small scales in the gluing theorem.

Joyce's theorem gives a unique Ricci-flat QALE metric $g$ with
K\"ahler form $\omega$ in each chosen class. On the curve model this
is the class of the averaged form. A holomorphic transition pulls $g$
back to a Ricci-flat QALE metric in the same class, so it preserves $g$
by uniqueness. For an antiholomorphic transition $\gamma$, the
K\"ahler form of $\gamma^*g$ is $-\gamma^*\omega$, again in the
chosen class. Uniqueness then gives
\[
 \gamma^*g=g,\qquad \gamma^*\omega=-\omega.
\]
In the fixed local product identifications, the limiting K\"ahler
classes are restrictions of the chosen class. This follows by
passing to the limit in the periods over the exceptional spheres.
Since $H^2(\EH;\R)\cong\R$, these periods determine the classes.
The normal complex structures are fixed, so uniqueness identifies
the Ricci-flat limits with the chosen Eguchi--Hanson metrics
\cite[Proposition~7.1 and the proof of Theorem~3.3]{QALE}.

On each complex normal resolution $\pi$ of dimension $m$, crepancy gives
a nowhere-vanishing holomorphic extension $\Omega$ of $\pi^*\Omega_{\rm flat}$.
Joyce's equation $\omega^m=C_m\,\Omega\wedge\bar\Omega$ has the same
constant $C_m$ as in the flat model~\cite[pp.~129--130]{QALE}.
Thus $\Omega$ has constant norm and is parallel. Its product with the flat holomorphic volume form on the tangent factor transforms
as in Lemma~\ref{lem:phase}. The transitions therefore preserve the associated
Cayley forms.

In each holomorphic comparison chart $\psi$, let $\omega,\Phi$ and
$\omega_0,\Phi_0$ denote the K\"ahler and Cayley forms on the resolution
and the product model, respectively, including the flat tangent factor
in both. The holomorphic volume forms agree under $\psi$, since they
extend the same flat form. Thus~\eqref{eq:phase} gives
\[
 \psi^*\Phi-\Phi_0=\omega_0\wedge D+\tfrac12 D\wedge D,
 \qquad D=\psi^*\omega-\omega_0.
\]
The product connection preserves the complex structure and $\omega_0$,
so the metric estimate~\cite[equation~(9)]{QALE} bounds every derivative
of $D$ with the same weights. The Leibniz rule and boundedness of the
zeroth-order weight give the required estimate for the Cayley form
\cite[Definition~13.1.8, equation~(13.3), p.~344]{JoyceBook}.

The curve resolution is used with its flat complex tangent line,
and the point resolution is $\EH\times\EH$.
Their contraction maps are the local morphisms of
Lemma~\ref{lem:markedresolution}. For every proper inclusion in the
singular part of $\mathcal L$ described in Section~\ref{sec:stabilizers},
the matched periods and uniqueness above give the required $A_1$
product limit.
The member $T^8\in\mathcal L$ carries the trivial normal resolution
and the original flat $\Spin(7)$ structure; compatibility with this
member is the flat-limit condition in the QALE asymptotics.

\begin{proposition}[Existence for the specified R-data]\label{prop:existence}
The complex models and the untwisted real model above define R-data
for the original flat quotient $\Omm$. The manifold $\Mmm$ constructed
in Section~\ref{sec:realL} admits a torsion-free $\Spin(7)$ structure.
\end{proposition}
\begin{proof}
Sections~\ref{sec:fan}--\ref{sec:globalV} provide the complex local
resolutions, and the preceding metric choices give compatible
Ricci-flat metrics and equivariant Cayley forms. The disjoint real
component has the untwisted R-data of Section~\ref{sec:realL}.
We use the local product identifications and lifted affine
transitions constructed above in Joyce's construction
\cite[Definition~13.2.5, p.~348]{JoyceBook}.
The resulting resolution is the manifold $\Mmm$ of
Section~\ref{sec:realL}: the local identifications commute with the
projections to $\Omm$, so on overlaps they agree over the dense
regular locus and hence everywhere by continuity.
Joyce's Theorem~13.6.2~\cite[p.~362]{JoyceBook} therefore applies
to the original flat quotient $\Omm=T^8/\Gamma$.
\end{proof}

\section{Factorization of the resolution map}\label{sec:comparison}
We construct the lift $\widetilde\tau$ on $Y$ and compute its action
on the exceptional classes used in Proposition~\ref{prop:abstracttransfer}.
We then show that $\pi_M$ factors through the quotient of $Y\times S^1$ by
$(\widetilde\tau,-1)$.

\subsection{The seven-dimensional resolution}
Let $\rho:Y\to X=T^7/G$ be the flat Kummer realization of the
resolution in \cite[Theorem~3 and Remark~4]{MM}. We use the standard
blow-down maps in the flat tubular coordinates. For each of the ten old singular components $N_j$ of $X$, the table
gives a representative torus in $T^7$. Its product with the $s$-circle
is a component of the fixed set in $T^8$ of the indicated involution
from~\eqref{eq:OUJ}. Parameters $a,b,c$ have period one, and $r=0,1,2,3$ corresponds in order to $(\epsilon_1,\epsilon_2)=(0,0),(0,\frac12),(\frac12,0),(\frac12,\frac12)$:
\begin{center}\small
\begin{tabular}{@{}clclc@{}}\toprule
Center & Fixing involution & $t$ & $(z_1,z_2,z_3)$ & Pullback of $\varphi$\\\midrule
$N_1$ & $R_{1,0}=F^3\kappa$ & $0$ & $(a-ia,b-ib,ic)$ & $2\,da\,db\,dc$\\
$N_2$ & $R_{1,1}=F^3\kappa\iota_1$ & $0$ & $(a+ia,b+ib,-\tfrac14+ic)$ & $-2\,da\,db\,dc$\\
$N_{3+r}$ & $O_{0,0}=\kappa$ & $\tfrac12$ & $(a+i\epsilon_1,b+i\epsilon_2,c)$ & $da\,db\,dc$\\
$N_{7+r}$ & $O_{0,1}=\kappa\iota_2$ & $\tfrac12$ & $(\epsilon_1+ia,\epsilon_2+ib,c-\tfrac i4)$ & $-da\,db\,dc$\\\bottomrule
\end{tabular}
\end{center}
Each torus in the table is a double cover of its image $N_j$. Its deck
involution is $\iota(a,b,c)=(-a,-b,c+\tfrac12)$, induced by $\iota_2$
for $j=1,2$ and by $\iota_1$ for $3\leq j\leq10$. With each $N_j$
oriented by $\varphi$, the last column gives
\begin{equation}\label{eq:Nvolumes}
 \int_{N_1}\varphi=\int_{N_2}\varphi=1,
 \qquad\int_{N_j}\varphi=\tfrac12\quad(3\leq j\leq10).
\end{equation}
Pullback identifies $H^*(N_j;\R)$ with the invariant cohomology of the
covering torus. Since $da$ and $db$ change sign while $dc$ is fixed,
a basis is $1$, $[dc]$, $[da\wedge db]$,
and $[da\wedge db\wedge dc]$. Thus $b(N_j)=(1,1,1,1)$.
The form $dc$ descends to $\ell_j$, which spans the space of parallel
one-forms on $N_j$. In the coordinates of the table,
\begin{equation}\label{eq:ell}
 \ell_j=dy_3\quad(j=1,2),\qquad
 \ell_j=dx_3\quad(3\leq j\leq10).
\end{equation}
Let $v_j=\ell_j^\sharp$ be the metric dual of $\ell_j$. The normal
complex structures $I_j$ are determined by
$g(I_j u,w)=(\iota_{v_j}\varphi)(u,w)=\varphi(v_j,u,w)$ for vectors
$u,w$ normal to $N_j$, as in~\cite[Section~2.2]{MM}.

The replacement near $N_j$ is $(\widehat N_j\times\EH)/\langle\iota_*\rangle$,
where $\widehat N_j$ is the torus cover. The differential of the ambient
involution acts by $+1$ on the $I_j$-complex line spanned by $\partial_t$
and by $-1$ on the normal two-plane in the $z_1,z_2$ directions.
In a parallel unitary frame adapted to these lines, the normal action
and the restriction of its lift to the exceptional sphere are
\[
 \iota_*(p,w_1,w_2)=(\iota p,w_1,-w_2),\qquad
 \iota_*(p,[W_1:W_2])=(\iota p,[W_1:-W_2]).
\]
The exceptional set $E_j=\rho^{-1}(N_j)\subset Y$ is a trivial oriented
sphere bundle over $N_j$. To check this explicitly, let $c$ be the circle coordinate translated by $1/2$. The map
\begin{equation}\label{eq:spheretrivial}
 [(p,[W_1:W_2])]\longmapsto
 ([p],[W_1:e^{2\pi i c(p)}W_2])
\end{equation}
is well defined, since the two minus signs cancel under $\iota$. Its inverse uses the opposite exponential. Orient $E_j$ by the given orientation of $N_j$ and the complex orientation
of the fibers. Write $\xi_j=\PD_Y[E_j]\in H^2(Y)$ for its Poincar\'e dual.

The flat Kummer construction and the resolution in~\cite[Section~2.2,
pp.~1115--1116]{MM} use the same pieces over each $N_j$. Their
attaching maps differ by a positive fiber rescaling and an equivariant
flow, hence are isotopic. Extending the isotopy across the collars
gives an orientation-preserving diffeomorphism that preserves each
$E_j$ and its coorientation. Thus $\xi_j$ corresponds to the
published class $x_j$.

\subsection{Lifting \texorpdfstring{$\tau$}{tau} and determining root signs}\label{sec:taulift}
The affine involution $\tau$ defined in~\eqref{eq:tau} normalizes $G$ by~\eqref{eq:relations}, so it descends to $X=T^7/G$. To compute its action on the specified representatives, compose it with $\kappa\iota_1\in G$. Direct substitution gives
\begin{center}\small
\begin{tabular}{@{}cccc@{}}\toprule
Source and target & Parameter map & Sign of $\ell_j$ & Sign of root\\\midrule
$N_1\to N_2$ & $(a,b,c)\mapsto(-a,-b,c)$ & $+$ & $-$\\
$N_2\to N_1$ & $(a,b,c)\mapsto(-a,-b,c)$ & $+$ & $-$\\
$N_{3+r}\to N_{3+r}$ & $(a,b,c)\mapsto(-a,-b,\tfrac34-c)$ & $-$ & $+$\\
$N_{7+r}\to N_{7+r}$ & $(a,b,c)\mapsto(a,b,\tfrac34-c)$ & $-$ & $+$\\\bottomrule
\end{tabular}
\end{center}
Equalities are modulo the torus periods; in the last row $-\epsilon_i=\epsilon_i$ modulo one. The table describes the induced maps of the quotient centers.

Because $\tau^*\varphi=-\varphi$, the transformation sign of the normal K\"ahler form is the negative of the transformation sign of $\ell_j$. Thus on the first pair the normal map is antiunitary and reverses the complex orientation of the exceptional $\CP^1$; on the remaining centers it is unitary and preserves that orientation. The minimal $A_1$ resolution is functorial for unitary transformations and for unitary transformations composed with conjugation. Its standard lifts commute with the radial identification outside the exceptional sphere and preserve the Eguchi--Hanson metric when the sizes agree. Choose equal sizes on the exchanged pair. These lifts and the original affine map outside the tubular neighborhoods give a smooth map $\widetilde\tau:Y\to Y$.

\begin{lemma}\label{lem:tauaction}
The map $\widetilde\tau$ is an orientation-reversing involution and
\begin{equation}\label{eq:tauroots}
 \widetilde\tau^*\xi_1=-\xi_2,\quad
 \widetilde\tau^*\xi_2=-\xi_1,\quad
 \widetilde\tau^*\xi_j=\xi_j\quad(3\leq j\leq10).
\end{equation}
\end{lemma}
\begin{proof}
The root signs are those of the preceding table. The lift squares to the identity on the complement of the exceptional sets, since it covers $\tau^2=1$. Its square is the identity everywhere by continuity; the same argument verifies equivariance with the normalizer action used to define each local bundle. Orientation reversal holds on the dense old smooth locus and therefore on $Y$.
\end{proof}
Only the smooth involution and its action on cohomology are needed below; no invariance of the perturbed $G_2$ metric is required.

For $j\geq3$, the fixing involution lies in $H$, so the corresponding
branch is resolved in $V$. Here $\ell_j=dx_3$, and contracting~\eqref{eq:phi} gives
\begin{equation}\label{eq:normalmatch}
 \iota_{\partial_{x_3}}\varphi
 =-dt\,dy_3+dx_1dx_2-dy_1dy_2.
\end{equation}
On the normal four-plane of each chosen representative, $ds\,dx_3$ vanishes,
so the displayed form is the restriction of $\omega'$. Thus $I_j$ agrees
there with the normal complex structure from $V$. On an $F$-translate
the two structures are opposite, since $F$ transports $I_j$ and
$F^*\omega'=-\omega'$ by Lemma~\ref{lem:phase}. Opposite complex
structures have the same complex lines $\Lambda$ in the real normal
space $W$. The identity on pairs $([w],\Lambda)$ with $w\in\Lambda$
therefore identifies their real $A_1$ resolutions over $W/\{\pm1\}$.
It commutes with the transition maps and reverses the complex
orientation of the exceptional sphere.

The tori $\widehat N_1\times S^1$ and $\widehat N_2\times S^1$ have the
same image in the smooth locus of $T^8/H$. In $V/\langle\bar F\rangle$
this becomes the real singular component $L$ of~\eqref{eq:L}, resolved
using $\ell=dy_3$.

\subsection{The intermediate quotient}\label{sec:localcomparison}
Define the compact quotient
\begin{equation}\label{eq:Z}
 Z=(Y\times S^1)/\langle\widetilde\sigma\rangle,
 \qquad \widetilde\sigma(y,s)=(\widetilde\tau y,-s).
\end{equation}
The product map $\rho\times\id\colon Y\times S^1\to X\times S^1$
resolves all the old branches. Since
$\rho\circ\widetilde\tau=\tau\circ\rho$, it descends to the partial
resolution $\bar\rho\colon Z\to\Omm$, given by
$\bar\rho([y,s])=[\rho(y),s]$.

\begin{proposition}[The degree-one map]\label{prop:degreeone}
For the R-data of Proposition~\ref{prop:existence}, the resolution map
$\pi_M$ of Section~\ref{sec:realL} factors as
\begin{equation}\label{eq:rmap}
 \Mmm\xrightarrow{\ r\ }Z\xrightarrow{\ \bar\rho\ }\Omm,
\end{equation}
where $r$ is a continuous resolution map of degree one. The singular
locus of $Z$ is a disjoint union of smooth four-dimensional centers
with transverse model $A_1=\C^2/\{\pm1\}$.
\end{proposition}
\begin{proof}
For both resolutions we use the flat tubular charts inherited from $T^8$
and the same radial identifications on their common collars.
The identification after~\eqref{eq:normalmatch} commutes with the
residual involution on the resolved old branch. At each intersection
the local model of $Z$ is therefore the first step of
Lemma~\ref{lem:markedresolution}. In the curve model
we contract the new exceptional spheres over the two disjoint centers.
In the point model the map is $\id\times\rho_{\EH}\colon
\EH\times\EH\to\EH\times A_1$, with the old factor first.
Away from the intersections we use $\rho_{\EH}$ on each ordinary new
branch and the same real identification on each ordinary old complex branch.

On the disjoint old odd-$F$ pair, $\tau$ exchanges $N_1$ and $N_2$
in $X$, so $\widetilde\tau$ exchanges their resolved neighborhoods in
$Y$. The quotient identifies the two product neighborhoods and is
smooth there. With the normal direction $\ell=dy_3$ and the sphere
size fixed in Section~\ref{sec:realL}, this is the same untwisted
replacement along $L$.

The new centers and their contractions are intrinsic to the marked
models, so the local maps commute with the transitions of
Section~\ref{sec:globalV}. On each overlap they agree over the dense
preimage of the regular locus of $\Omm$, hence everywhere by continuity.
They glue to a continuous map $r$ with $\bar\rho\circ r=\pi_M$.
The local models give the stated centers and transverse singularities.
Over the regular locus of $\Omm$ the map $r$ is an orientation-preserving
diffeomorphism. Since $Z$ has a rational fundamental class,
$r_*[\Mmm]=[Z]$.
\end{proof}

\section{Simple connectivity, nonformality and holonomy}\label{sec:topology}
\subsection{Simple connectivity}
\begin{proposition}\label{prop:pi1M}
The manifold $\Mmm$ is simply connected.
\end{proposition}
\begin{proof}
The manifold $Y$ is simply connected by
\cite[Theorem~3, pp.~1113--1114]{MM}.
The involution $\widetilde\tau$ has a fixed point in the old smooth locus:
take $t=1/2$, $z_3=1/8$, and generic $z_1,z_2$.
On the simply connected space $Y\times\R$, let
\[
 T(y,s)=(y,s+1),\qquad R(y,s)=(\widetilde\tau y,-s).
\]
These transformations generate the infinite dihedral group, with
$R^2=1$ and $RTR=T^{-1}$, and the orbit space is $Z$.
Both $R$ and $TR$ have fixed points, at $s=0$ and $s=1/2$ over the
chosen point of $Y$. Armstrong's orbit-space theorem~\cite{Armstrong}
identifies $\pi_1(Z)$ with this group modulo the normal subgroup
generated by elements with fixed points. Since $R$ and $TR$ generate
the group, $\pi_1(Z)=1$.

By Proposition~\ref{prop:degreeone}, $\Mmm\to Z$ replaces disjoint
normal $A_1$ cones by Eguchi--Hanson spaces. The common boundary has
fiber $\RP^3$. Both fillings have simply connected fibers, so their
neighborhoods have the fundamental group of the same base center,
and the boundary maps induce the same homomorphism to that group.
Van Kampen therefore gives $\pi_1(\Mmm)=\pi_1(Z)=1$, including in
the presence of flat bundle monodromy.
\end{proof}

\subsection{Nonformality}\label{sec:nonformality}
For the exceptional classes $\xi_j=\PD_Y[E_j]$, set
\begin{equation}\label{eq:abcw}
 a=\xi_1+\xi_2,\qquad b=\xi_7+\xi_3,\qquad
 c=\xi_7-\xi_3,\qquad w=\xi_1-\xi_2.
\end{equation}
The cohomology relations and the primitive constructed by
\MM\ give
\begin{equation}\label{eq:fourzero}
 ab=bc=aw=cw=0,\qquad
 \langle mw,[Y]\rangle=16
 \quad(m\in\langle a,b,c\rangle)
\end{equation}
in our sign convention; see \cite[Lemma~16 and Theorem~17]{MM}.
Appendix~\ref{sec:mm} records the normalization and checks this period.
By Lemma~\ref{lem:tauaction},
\[
 \widetilde\tau^*(a,b,c,w)=(-a,b,c,w).
\]
Together with the degree-one map of Proposition~\ref{prop:degreeone},
these identities satisfy the hypotheses of
Proposition~\ref{prop:abstracttransfer}.
The descended classes $a\,ds,b,c,w$, pulled back to $\Mmm$,
therefore give~\eqref{eq:mainperiod}, with pairing $16/2=8$ for
every representative. Thus $\Mmm$ is nonformal over both $\Q$ and $\R$.

\subsection{Full holonomy}
\begin{proposition}\label{prop:fullholonomy}
The torsion-free metric on $\Mmm$ has holonomy $\Spin(7)$.
\end{proposition}
\begin{proof}
The metric has holonomy contained in $\Spin(7)$ by
Proposition~\ref{prop:existence}, and $\Mmm$ is simply connected by
Proposition~\ref{prop:pi1M}. Joyce's classification
\cite[Theorem~10.6.1, p.~259]{JoyceBook} leaves four possibilities:
\[
 \Spin(7),\qquad \mathrm{SU}(4),\qquad
 \mathrm{Sp}(2),\qquad \mathrm{SU}(2)\times\mathrm{SU}(2).
\]
In each of the last three cases, a parallel complex structure would
make $\Mmm$ K\"ahler and hence formal over $\R$ by~\cite{DGMS}.
This contradicts the Massey product above. Thus the holonomy is
$\Spin(7)$.
\end{proof}
This completes the proof of Theorem~\ref{thm:main}.

\appendix
\section{The seven-dimensional Massey period}\label{sec:mm}
We use the published calculation of~\cite{MM} and record its
normalization in the conventions of this paper.
Let $\eta_j$ denote the Thom $2$-form of $E_j\subset Y$ constructed
in~\cite[equation~(18)]{MM}. These forms have pairwise disjoint
supports and satisfy $[\eta_j]=\xi_j$. On each exceptional fiber
$S\cong\CP^1$, the normal bundle of $E_j$ in $Y$ restricts to
$\mathcal O_S(-2)$. Thus $\xi_j|_{\EH}=-2u$ for the generator
of Section~\ref{sec:ehmarkings}, and
\begin{equation}\label{eq:Thomnorm}
 \int_{\CP^1}\eta_j=-2,\qquad
 \int_{\EH}\eta_j^2=-2.
\end{equation}
The cohomology relations in~\cite[Lemmas~9--10, equations~(20)--(21)]{MM}
are
\begin{equation}\label{eq:MMrelations}
 \xi_i\xi_j=0\ (i\ne j),\qquad \xi_1^2=\xi_2^2,
 \qquad \xi_3^2=\xi_4^2=\cdots=\xi_{10}^2.
\end{equation}
They give $ab=bc=aw=cw=0$ for~\eqref{eq:abcw}.
For the disjoint-support representatives, $ab=0$ as a form and
$bc=\eta_7^2-\eta_3^2$.

By~\cite[Lemma~16, p.~1125]{MM}, there is a smooth three-form $v$ with
\begin{equation}\label{eq:vlocal}
 \begin{aligned}
 dv&=\eta_7^2-\eta_3^2,\\
 v\eta_1&=-4\pi_1^*(\varphi|_{N_1})\eta_1,\qquad
 v\eta_2=4\pi_2^*(\varphi|_{N_2})\eta_2,
 \end{aligned}
\end{equation}
where $\pi_j$ is the projection of the resolved neighborhood onto
$N_j$. Taking the primitive of $ab$ to be zero, our
convention~\eqref{eq:masseyconvention} gives $m=av$.
If $e_j=[\pi_j^*(\varphi|_{N_j})\eta_j]$, then
$[m]=-4e_1+4e_2$. Disjoint supports and fiber integration give
\[
 \langle e_j\xi_k,[Y]\rangle=0\quad(j\ne k),\qquad
 \langle e_j\xi_j,[Y]\rangle=-2\int_{N_j}\varphi.
\]
Using~\eqref{eq:Nvolumes}, we obtain
\begin{equation}\label{eq:period16}
 \langle mw,[Y]\rangle=(-4)(-2)-4(-2)=16.
\end{equation}
Since $aw=cw=0$, Lemma~\ref{lem:detector} gives this value for every
real representative. The convention of~\cite[equation~(23) and
Theorem~17]{MM} uses $-m$ and gives $-16$.

The classes $\xi_j$ are integral, so injectivity of
$H^*(Y;\Q)\to H^*(Y;\R)$ gives $ab=bc=aw=cw=0$ over $\Q$ as well.
Thus rational defining systems exist, and each extends to a real one
whose pairing is $16$ by~\eqref{eq:period16}.

\section{Betti numbers and signature}\label{sec:manualorbits}
We prove Corollary~\ref{cor:numerical} by computing component orbits and their tangent representations. The procedure is uniform: label the fixed components by binary coordinates, compute the induced translations, and take invariant exterior monomials for the stabilizer of one component. The final step uses full holonomy to determine the signature.

\subsection{The elementary abelian subgroup and conjugation}\label{sec:manualH}
Write an element of $H\cong(\Z/2)^5$ as
\begin{equation}\label{eq:Hbits}
 v=F^{2p}\kappa^b\iota_1^r\iota_2^\ell\sigma^e,
 \qquad p,b,r,\ell,e\in\mathbb F_2.
\end{equation}
From~\eqref{eq:Aaction}--\eqref{eq:Baction}, its action is
\begin{align}
 z_j&\longmapsto(-1)^{p+r+\ell}\operatorname{conj}^{b}(z_j),\quad j=1,2,\nonumber\\
 t&\longmapsto(-1)^{b+e}\bigl(t-\tfrac12\bigr)+\tfrac12+2p,\nonumber\\
 x_3&\longmapsto(-1)^e x_3+e/4+r/2,\label{eq:Haction}\\
 y_3&\longmapsto(-1)^{b+e}y_3+(-1)^b\ell/2,\nonumber\\
 s&\longmapsto(-1)^e s.\nonumber
\end{align}
The elements of $H$ with fixed points are the families $O,U,J$ of~\eqref{eq:OUJ}.

\begin{lemma}\label{lem:Fpairing}
Conjugation by $F$ acts on these families by
\begin{equation}\label{eq:Fconjugation}
 O_{\alpha d}\longmapsto O_{\alpha+1,d},\qquad
 U_{\alpha cd}\longmapsto U_{\alpha+1,c+1,d},\qquad
 J_c\longmapsto J_{c+1}.
\end{equation}
Consequently every complex exceptional sector is exchanged with another by $\bar F$. No odd-$F$ element preserves a maximal Klein stabilizer or a fixed component belonging to the families $O,U,J$.
\end{lemma}
\begin{proof}
The generator identities are
\[
 F\kappa F^{-1}=F^2\kappa,\qquad
 F\sigma F^{-1}=F^2\iota_1\sigma,\qquad
 F\iota_jF^{-1}=\iota_j.
\]
They give~\eqref{eq:Fconjugation}. Since $H$ is abelian, conjugation by every element of $FH$ is the same as conjugation by $F$ on $H$. Each of the elements $O_{\alpha d}$, $U_{\alpha cd}$ and $J_c$ is changed by conjugation. At a surface component its generic stabilizer is generated by that element, so a setwise stabilizer must preserve the element. At a maximal Klein stabilizer the unique old element is $O_{\alpha d}$, which is also changed by~\eqref{eq:Fconjugation}. The same argument applies to the isolated fixed loci, whose full stabilizers are maximal Klein groups.
\end{proof}

\subsection{Component labels and their translations}\label{sec:labeltranslations}
For an old element $O_{\alpha d}$, the fixed components are described by
\begin{equation}\label{eq:Olabels}
 t=\alpha+\tfrac12+2T,\qquad
 y_3=-\tfrac d4+\tfrac Y2,\qquad T,Y\in\mathbb F_2.
\end{equation}
If $\alpha+d=0$, the two imaginary coordinates $y_1,y_2$ lie in $\{0,1/2\}$ and $x_1,x_2$ are free; if $\alpha+d=1$, the two real coordinates lie in $\{0,1/2\}$ and the imaginary coordinates are free. Let $\epsilon\in\mathbb F_2^2$ record those two fixed coordinates. The remaining coordinates $x_3,s$ are free. Thus there are $2^4=16$ components.

For $U_{\alpha cd}$, use labels $T,X,Y,S\in\mathbb F_2$ for the fixed $\mathsf B$ coordinates:
\begin{equation}\label{eq:Ulabels}
 t=\alpha+\tfrac12+2T,\quad
 x_3=\tfrac18+\tfrac c4+\tfrac X2,\quad
 y_3=\tfrac d4+\tfrac Y2,\quad s=\tfrac S2.
\end{equation}
When $\alpha+c+d=0$, all four $\mathsf A$ coordinates are free and there are $16$ surface components. When $\alpha+c+d=1$, the four $\mathsf A$ coordinates lie in $\{0,1/2\}$, giving $16\cdot16=256$ points.

For $J_c$, the fixed coordinates in $\mathsf B$ are
\begin{equation}\label{eq:Jlabels}
 x_3=\tfrac18+\tfrac c4+\tfrac X2,\qquad s=\tfrac S2,
\end{equation}
while $t,y_3$ are free. Its two fixed $\mathsf A$ coordinates are the imaginary coordinates for $c=0$ and the real coordinates for $c=1$. They give two further bits $\epsilon\in\mathbb F_2^2$, hence $16$ surface components.

\begin{samepage}
\begin{lemma}[Label translations]\label{lem:labeltranslations}
For the element $v$ in~\eqref{eq:Hbits}, all fixed $\mathsf A$ labels and the bit $S$ are unchanged. The remaining labels change by the following translations in $\mathbb F_2$:
\begin{center}
\begin{tabular}{@{}clc@{}}\toprule
Family & Translation vector & Rank\\\midrule
$O_{\alpha d}$ & $(\Delta T,\Delta Y)=(p+\alpha(b+e),\ \ell+d(b+e))$ & $2$\\
$U_{\alpha cd}$ & $(\Delta T,\Delta X,\Delta Y)=(p+\alpha(b+e),\ r+ce,\ \ell+d(b+e))$ & $3$\\
$J_c$ & $\Delta X=r+ce$ & $1$\\\bottomrule
\end{tabular}
\end{center}
The translations are independent of the initial component labels.
\end{lemma}
\end{samepage}
\begin{proof}
Every fixed $\mathsf A$ coordinate is $0$ or $1/2$, and the action of $H$ is diagonal by signs, so those labels are unchanged. Also $-S/2=S/2$ modulo one. For the $t$ coordinate in~\eqref{eq:Olabels} or~\eqref{eq:Ulabels}, substitution into~\eqref{eq:Haction}, followed by division of the difference by two, gives
$\Delta T=p+\alpha(b+e)$ modulo two. For example, substitute $x_3=1/8+c/4+X/2$ from the family $U$ into~\eqref{eq:Haction}. Writing the image in the same form gives
\[
 X'\equiv\begin{cases}X+r,&e=0,\\-X+r-c,&e=1\end{cases}\pmod2,
 \qquad X'-X=r+ce\quad\text{in }\mathbb F_2.
\]
This also proves the formula for $J_c$. Substitution of $y_3$ gives $\Delta Y=\ell+d(b+e)$ for both choices of the constant $\pm d/4$. These expressions do not contain the starting bits. Finally, varying $p$, $r$, and $\ell$ separately gives the independent coordinate translations in the rows where they occur. This proves the stated ranks.
\end{proof}
A translation action of rank $r$ has orbits of size $2^r$ and kernel of order $2^{5-r}$. Thus the component stabilizers are
\begin{align}
 N_O&=\{p=\alpha(b+e),\ \ell=d(b+e)\},\label{eq:Ostabilizer}\\
 N_U&=\{p=\alpha(b+e),\ r=ce,\ \ell=d(b+e)\},\label{eq:Ustabilizer}\\
 N_J&=\{r=ce\}.\label{eq:Jstabilizer}
\end{align}
These equations give the stabilizer for every component in the family.

\subsection{The cohomology of normalized surface closures}\label{sec:tangentcharacters}
For each family, the tangent representation follows by substituting the stabilizer equations into~\eqref{eq:Haction}. For $O_{\alpha d}$, put $\eta=\alpha+d$. On the two free $\mathsf A$ coordinates the common sign is $(-1)^{r+\eta e}$, whereas on $x_3,s$ it is $(-1)^e$. The two characters $r+\eta e$ and $e$ are independent on $N_O$. The invariant exterior monomials therefore contain either both or neither of the directions in each pair. Its Poincar\'e polynomial is $(1+q^2)^2$.

For a surface $U_{\alpha cd}$, the equation $\alpha+c+d=0$ gives
$p+r+\ell=cb$ on $N_U$. Its tangent action on $(z_1,z_2)$ is consequently $(-1)^{cb}\operatorname{conj}^b$. This representation has two trivial real directions and two directions with the same nontrivial sign character. Its invariant exterior polynomial is
\[
 (1+q)^2(1+q^2)=1+2q+2q^2+2q^3+q^4.
\]
For $J_c$, the two free $\mathsf A$ directions have sign $(-1)^{p+\ell+c(b+e)}$ and the $t,y_3$ pair has sign $(-1)^{b+e}$. These characters are independent on $N_J$, giving again $(1+q^2)^2$.

Translations in these component stabilizers act trivially on torus cohomology. Rational cohomology of a finite quotient is the invariant cohomology. Therefore these tangent calculations give the cohomology of every normalized surface closure, including those whose stabilizer action has fixed points. The resulting orbit counts and invariant polynomials are
\begin{center}\small
\begin{tabular}{@{}lccccl@{}}\toprule
Family & Elements & Components & Orbit & Orbits/element & Polynomial\\\midrule
$O_{\alpha d}$ & $4$ & $16$ & $4$ & $4$ & $(1+q^2)^2$\\
$U_{\alpha cd}$, even & $4$ & $16$ & $8$ & $2$ & $(1+q)^2(1+q^2)$\\
$J_c$ & $2$ & $16$ & $2$ & $8$ & $(1+q^2)^2$\\
$U_{\alpha cd}$, odd & $4$ & $256$ & $8$ & $32$ & $1$ (point)\\\bottomrule
\end{tabular}
\end{center}
Here ``even'' and ``odd'' refer to $\alpha+c+d$. There are $16+16=32$ surface sectors of the first polynomial type, $8$ of the second, and $128$ point sectors. At a point exactly one nonidentity element has an isolated fixed locus, namely the odd $U$, so the last row counts point strata once rather than three times.

\subsection{The lower strata and the real component}
For $K(2\alpha,c,d)=\langle O_{\alpha d},U_{\alpha cd}\rangle$, the $\mathsf B$ levels are those of~\eqref{eq:Ulabels}. Any element of $H$ preserving a component must satisfy~\eqref{eq:Ustabilizer}. That kernel is exactly $K(2\alpha,c,d)$: its free parameters $b,e$ generate $O_{\alpha d}$ and $U_{\alpha cd}$. Each of its elements fixes the common fixed component pointwise. Thus its component normalizer in $H$ has order four. By Lemma~\ref{lem:Fpairing} there is no further normalizer in $\Gamma$. Each $\Gamma$ orbit has size $16$. The four even triples give $4\cdot64/16=16$ curve strata, and the four odd triples give $4\cdot256/16=64$ point strata.

The old odd-$F$ components form one orbit, by the normalizer calculation
in Section~\ref{sec:realL}. The $40$ complex surface orbits under $H$
pair under $F$, so there are $40/2+1=21$ open four-dimensional strata.
Removing the lower-dimensional intersections leaves these strata
connected, since the removed subsets have real codimension at least two.

\subsection{Additive cohomology of the complex resolution}\label{sec:cohomology}
An invariant real exterior monomial on $T^8$ contains an even number of directions from each of the four complex coordinate planes, or an odd number from each plane. These are the only parity vectors pairing trivially with all even sign changes. Consequently
\begin{equation}\label{eq:baseP}
 P_{T^8/H}(q)=(1+q^2)^4+16q^4
 =1+4q^2+22q^4+4q^6+q^8.
\end{equation}
The map $f:V\to T^8/H$ is semismall: fibers over complex codimensions $2,3,4$ have dimensions $1,1,2$, respectively. Hence the relevant strata are surfaces and points; curves contribute no separate summand. Normalized surface closures are finite holomorphic quotients of complex two-tori and are rational homology manifolds. Their intersection cohomology is therefore ordinary rational cohomology. The top-fiber local systems are trivial: holomorphic transitions preserve the fundamental class of $\CP^1$, and the point fiber $\CP^1\times\CP^1$ has one top-dimensional component. The finite normalization map supplies the intersection-complex extension if a surface closure has locally identified branches.

The projective semismall decomposition theorem~\cite[Theorem~3.4.1]{dCM} and the preceding four-row count give
\begin{align}
 P_V(q)&=P_{T^8/H}(q)+32q^2(1+q^2)^2\nonumber\\
 &\quad+8q^2(1+q)^2(1+q^2)+128q^4,\label{eq:PV}\\
 b(V)&=(1,0,44,16,230,16,44,0,1),\quad\chi(V)=288.\nonumber
\end{align}

\subsection{The real replacement}
For $L=N_1\times S^1$, the effective deck action on the torus cover
has tangent representation $\operatorname{diag}(-1,-1,1,1)$.
Its invariant exterior algebra gives
\begin{equation}\label{eq:Lbetti}
 P_L(q)=(1+q)^2(1+q^2)=1+2q+2q^2+2q^3+q^4.
\end{equation}
The root local system of the last replacement is trivial. Its global
exceptional class restricts to $-2$ on each exceptional sphere, so
Leray--Hirsch and Mayer--Vietoris across the $\RP^3$ boundary give
\begin{equation}\label{eq:lastcohom}
 H^k(\Mmm;\Q)\cong H^k(V/\bar F;\Q)\oplus H^{k-2}(L;\Q).
\end{equation}

\subsection{The residual action and the final numbers}
Every exceptional sector is paired by Lemma~\ref{lem:Fpairing}, so its contribution to the residual trace is zero. Indeed, in the semismall decomposition the summands are canonically grouped by their supports. The residual involution sends the group for one support to that for its image, acting also on the top homology of the fiber. On a pair of exchanged supports its action is off-diagonal and has trace zero. On the four untwisted degree-two area forms, $F$ exchanges the first two with minus signs and negates each of the last two. Thus its trace is $-2$. Its trace on $H^3(V)$ is zero. Since $\chi(\Fix(\bar F))=\chi(L)=0$, the Lefschetz formula gives
\[
 0=2+2(-2)-2(0)+\operatorname{tr}(\bar F|H^4(V)),
 \qquad\operatorname{tr}(\bar F|H^4(V))=2.
\]
Poincar\'e duality supplies the other traces. Taking invariants and then applying~\eqref{eq:lastcohom} yields
\begin{align}
 b(V/\bar F)&=(1,0,21,8,116,8,21,0,1),\nonumber\\
 b(\Mmm)&=(1,0,22,10,118,10,22,0,1).\label{eq:computedM}
\end{align}
The alternating sum gives $\chi(\Mmm)=144$.
Full $\Spin(7)$ holonomy implies $\widehat A(\Mmm)=1$
\cite[Theorem~10.6.1]{JoyceBook}. On an eight-manifold with a
$\Spin(7)$ structure,
\begin{equation}\label{eq:indexidentity}
 48\widehat A=3\sgn-\chi.
\end{equation}
Indeed, this follows by substituting
$p_1^2-4p_2+8e=0$ into the degree-eight formulas for $\widehat A$
and the $L$-genus. Therefore
\[
 \sgn(\Mmm)=\frac{144+48}{3}=64,\qquad
 b_4^+=\frac{118+64}{2}=91,\qquad b_4^-=27.
\]
This proves Corollary~\ref{cor:numerical}.

\section*{Acknowledgments}
I am grateful to Alexander Zuevsky, whose message prompted me to return to
the question of formality for manifolds with exceptional holonomy.


\begin{thebibliography}{99}
\bibitem{AK} M. Amann and V. Kapovitch,
\emph{On fibrations with formal elliptic fibers},
Advances in Mathematics \textbf{231} (2012), no.~3--4, 2048--2068.
\href{https://doi.org/10.1016/j.aim.2012.07.022}{doi:10.1016/j.aim.2012.07.022}.

\bibitem{AT} M. Amann and I.~A. Taimanov,
\emph{On the formality of nearly K\"ahler manifolds and of Joyce's examples in $G_2$-holonomy},
Moscow Mathematical Journal \textbf{24} (2024), no.~4, 495--512.
\href{https://doi.org/10.17323/1609-4514-2024-24-4-495-512}{doi:10.17323/1609-4514-2024-24-4-495-512}.

\bibitem{Armstrong} M. A. Armstrong,
\emph{The fundamental group of the orbit space of a discontinuous group},
Proceedings of the Cambridge Philosophical Society \textbf{64} (1968), 299--301.
\href{https://doi.org/10.1017/S0305004100042845}{doi:10.1017/S0305004100042845}.

\bibitem{Bazaikin2005} Ya.~V. Bazaikin,
\emph{Special K\"ahler metrics on complex line bundles and the geometry of $K3$-surfaces},
Siberian Mathematical Journal \textbf{46} (2005), no.~6, 995--1004.
\href{https://doi.org/10.1007/s11202-005-0095-x}{doi:10.1007/s11202-005-0095-x}.

\bibitem{Bredon} G. E. Bredon,
\emph{Introduction to Compact Transformation Groups},
Pure and Applied Mathematics, Vol.~46, Academic Press, New York--London, 1972.
\href{https://doi.org/10.1016/S0079-8169(08)60134-8}{Chapter~III: doi:10.1016/S0079-8169(08)60134-8}.

\bibitem{Cavalcanti} G. R. Cavalcanti,
\emph{Formality of $k$-connected spaces in $4k+3$ and $4k+4$ dimensions},
Mathematical Proceedings of the Cambridge Philosophical Society
\textbf{141} (2006), no.~1, 101--112.
\href{https://doi.org/10.1017/S0305004106009340}{doi:10.1017/S0305004106009340}.

\bibitem{dCM} M. A. A. de Cataldo and L. Migliorini,
\emph{The Hard Lefschetz Theorem and the topology of semismall maps},
Annales scientifiques de l'\'Ecole Normale Sup\'erieure (4)
\textbf{35} (2002), no.~5, 759--772.
\href{https://doi.org/10.1016/S0012-9593(02)01108-4}{doi:10.1016/S0012-9593(02)01108-4}.

\bibitem{DGMS} P. Deligne, P. Griffiths, J. Morgan and D. Sullivan,
\emph{Real homotopy theory of K\"ahler manifolds},
Inventiones Mathematicae \textbf{29} (1975), 245--274.
\href{https://doi.org/10.1007/BF01389853}{doi:10.1007/BF01389853}.

\bibitem{FKLS} D. Fiorenza, K. Kawai, H.~V. L\^e and L. Schwachh\"ofer,
\emph{Almost formality of manifolds of low dimension},
Annali della Scuola Normale Superiore di Pisa, Classe di Scienze (5)
\textbf{22} (2021), no.~1, 79--107.
\href{https://doi.org/10.2422/2036-2145.201905_002}{doi:10.2422/2036-2145.201905\_002}.

\bibitem{Illman} S. Illman,
\emph{Smooth equivariant triangulations of $G$-manifolds for $G$ a finite group},
Mathematische Annalen \textbf{233} (1978), no.~3, 199--220.
\href{https://doi.org/10.1007/BF01405351}{doi:10.1007/BF01405351}.

\bibitem{JoyceBook} D. D. Joyce,
\emph{Compact Manifolds with Special Holonomy},
Oxford Mathematical Monographs, Oxford University Press, 2000.
\href{https://doi.org/10.1093/oso/9780198506010.001.0001}{doi:10.1093/oso/9780198506010.001.0001}.

\bibitem{QALE} D. D. Joyce,
\emph{Quasi-ALE metrics with holonomy $\mathrm{SU}(m)$ and $\mathrm{Sp}(m)$},
Annals of Global Analysis and Geometry \textbf{19} (2001), no.~2, 103--132.
\href{https://doi.org/10.1023/A:1010778214851}{doi:10.1023/A:1010778214851}.

\bibitem{MM} L. \MM,
\emph{Compact holonomy $G_2$ manifolds need not be formal},
Geometry \& Topology \textbf{30} (2026), no.~3, 1109--1127.
\href{https://doi.org/10.2140/gt.2026.30.1109}{doi:10.2140/gt.2026.30.1109}.

\bibitem{MilivojevicSpin7}
\emph{A nonformal compact manifold with full holonomy $\Spin(7)$
obtained from the Mart\'in--Merch\'an $G_2$ construction}.
Posted by A.~Milivojevi\'c, September 2026.
\href{https://github.com/aleksmilivojevic/nonformal-holonomy-spin7/blob/43f997df2c64fc63dfc292b7831a699473886416/nonformal_spin7_construction.pdf}{GitHub manuscript (PDF)}.

\bibitem{Sullivan} D. Sullivan,
\emph{Infinitesimal computations in topology},
Publications Math\'ematiques de l'IH\'ES \textbf{47} (1977), 269--331.
\href{https://doi.org/10.1007/BF02684341}{doi:10.1007/BF02684341}.
\end{thebibliography}
\end{document}